%% file: paperRe.tex
\documentclass[journal, final, A4paper, onecolumn, 12pt]{IEEEtran}

\usepackage{colordvi}
\usepackage[dvipsnames]{xcolor}

\usepackage[cmex10]{amsmath}
\usepackage[T1]{fontenc}
\usepackage{letterspace}
\usepackage{amssymb}
\usepackage{amsfonts}
\usepackage{mathtools}
\usepackage{booktabs}
\usepackage{psfrag}
\usepackage{array}
\usepackage{multirow} 
\usepackage{fixltx2e}
\usepackage{cite}
\usepackage{graphicx}
\usepackage{commath}
\usepackage[font=footnotesize,caption=false]{subfig}
\usepackage{tikz, pgfplots}
\usetikzlibrary{shapes,fit,arrows,decorations.markings,matrix}
\usepackage[mediumspace,mediumqspace,squaren,binary]{SIunits}

\pgfplotsset{
	compat=newest, 
	/pgfplots/ylabel absolute/.style={%
		/pgfplots/every axis y label/.style={at={(0,0.5)},xshift=-15pt,rotate=90},
		/pgfplots/every y tick scale label/.style={at={(0,1)},above right,inner sep=0pt,yshift=0.3em}
	}
}%

\newcommand{\ie}{\emph{i.e.},}
\newcommand{\eg}{\emph{e.g.},}

\newcommand{\mc}{\mathcal}

\newcommand{\ts}[1]{{\textnormal{#1}}}

\title{Optimization-based Islanding of Power Networks using Piecewise Linear AC Power Flow}%
\author{P.~A.~Trodden,~\IEEEmembership{Member,~IEEE,} W.~A.~Bukhsh,~\IEEEmembership{Student Member,~IEEE,}
  A.~Grothey, and K.~I.~M.~McKinnon%
\thanks{Submitted to IEEE Transactions on Power Systems}
\thanks{This work was supported
  by the UK Engineering and Physical Sciences Research Council (EPSRC)
  under grant EP/G060169/1.}  
\thanks{P.~A.~Trodden is with the
  Department of Automatic Control \& Systems Engineering, University
  of Sheffield, Mappin Street, Sheffield S1 3JD, UK (e-mail:
  \texttt{p.trodden@shef.ac.uk}).}
\thanks{W.~A.~Bukhsh, A.~Grothey, and K.~I.~M.~McKinnon are with the
  School of Mathematics, University of Edinburgh, James Clerk Maxwell
  Building, Edinburgh EH9~3JZ, UK (e-mail:
  \texttt{w.a.bukhsh@sms.ed.ac.uk},
  \texttt{a.grothey@ed.ac.uk}, \texttt{k.mckinnon@ed.ac.uk}).}}%
\begin{document}
\maketitle

\begin{abstract}

In this paper, a flexible optimization-based framework for intentional
islanding is presented. The decision is made of which transmission
lines to switch in order to split the network while minimizing
disruption, the amount of load shed, or grouping coherent
generators. The approach uses a piecewise linear model of AC power
flow, which allows the voltage and reactive power to be considered
directly when designing the islands. Demonstrations on standard test
networks show that solution of the problem provides islands that are
balanced in real and reactive power, satisfy AC power flow laws, and
have a healthy voltage profile.
\end{abstract}


\section*{Nomenclature}

\subsection*{Sets}
\begin{description}[\IEEEsetlabelwidth{$P_g^{\ts{G}}, P_g^{\ts{G}+}$}\IEEEusemathlabelsep]
\item[$\mc{B}$] Buses.
\item[$\mc{L}$] Lines.
\item[$\mc{G}$] Generators.
\item[$\mc{D}$] Loads.
\item[$\mc{B}_l$] Buses connected by line $l$.
\item[$\mc{L}_i$] Lines connected to bus $i$.
\item[$\mc{G}_i$] Generators located at bus $i$.
\item[$\mc{D}_i$] Loads located at bus $i$.
\item[$\mc{B}^\ts{0}$] Buses assigned to section $0$.
\item[$\mc{B}^\ts{1}$] Buses assigned to section $1$.
\item[$\mc{L}^\ts{0}$] Set of uncertain lines.
\item[$\mc{B}^\ts{G}$] Set of generator buses.
\end{description}
\subsection*{Parameters}
\begin{description}[\IEEEsetlabelwidth{$P_g^{\ts{G}}, P_g^{\ts{G}+}$}\IEEEusemathlabelsep]
\item[$G^\ts{B}_i, B^\ts{B}_i$] Shunt conductance, susceptance at bus $i$.
\item[$g_l, b_l, b^\ts{C}_l$] Conductance, susceptance, shunt susceptance of line $l$.
\item[$\tau_l$] Off-nominal tap ratio of line $l$ (if transformer).
\item[$V_i^-, V_i^+$] Min., max. voltage magnitude at bus $i$.
\item[$P_g^{\ts{G}-}, P_g^{\ts{G}+}$] Min., max. real power outputs of generator $g$.
\item[$Q_g^{\ts{G}-}, Q_g^{\ts{G}+}$] Min., max. reactive power outputs of generator $g$.
\item[$P_d^\ts{D}, Q_d^\ts{D}$] Real, reactive power demands of load $d$.
\item[$P^{\ts{L}+}_l$] Real power loss limit of line $l$.
\item[$\Theta_l, \Theta^+_l$] Max. angle across $l$ if connected, disconnected.
\item[$c_g(p^\ts{G}_g)$] Generation cost function for generator $g$.
\item[$\beta_d$] Loss penalty for load $d$.

\end{description}
\subsection*{Variables}
\begin{description}[\IEEEsetlabelwidth{$p_l^{ij}, q_l^{ij}$}\IEEEusemathlabelsep]
\item[$v_i, \delta_i$] Voltage magnitude and phase at bus $i$.
\item[$\theta_{ij}$] $\delta_{i} - \delta_{j}$, voltage phase
  difference between bus $i$ and $j$. Note $\theta_{ij} =
  -\theta_{ji}$.
\item[$y_{ij}$] $\cos \theta_{ij}$. Note $y_{ij} = y_{ji}$.
\item[$z_{ij}$] $\sin \theta_{ij}$. Note $z_{ij} = -z_{ij}$.
\item[$v^i_l, v^j_l$] Voltage magnitudes at either end of line $l$ (which
  connects buses $i$ and $j$).
\item[$\theta^{ij}_l$] Voltage phase difference across a line
  $l$. Note $\theta^{ij}_l = -\theta^{ji}_l$.
\item[$y^{ij}_l$] $\cos \theta^{ij}_l$. Note $y^{ij}_l = y^{ji}_l$.
\item[$p_l^{ij}, q_l^{ij}$] Real, reactive power injection at bus $i$
  into line $l$ (which connects buses $i$ and $j$).
\item[$p_g^\ts{G}, q_g^\ts{G}$] Real, reactive power outputs of generator $g$.
\item[$p_d^\ts{D}, q_d^\ts{D}$] Real, reactive power supplied to load $d$.
\item[$\alpha_d$] Proportion of load $d$ supplied.
\item[$\gamma_i$] Binary. Section ($0$ or $1$) assignment of bus $i$.
\item[$\zeta_g$] Binary. Connection status of generator $g$.
\item[$\rho_l$] Binary. Connection status of line $l$.
\end{description}

\section{Introduction}

\IEEEPARstart{T}{he} last decade has seen a number of notable cases of
wide-area blackouts as a consequence of severe disturbances and
cascading failures~\cite{ItalyUCTEreport,USblackout04,LD06}. Although
preventive and corrective systems exist to ameliorate the effects of
severe disturbances, the operation of networks closer to limits,
together with increased uncertainty in load and distributed
generation, means that cascading failures may be harder to prevent, or
stop once instigated~\cite{Bialek05}. Thus, intentional islanding is
attracting attention as a corrective measure for limiting
the effects of severe disturbances and preventing wide-area blackout.

Intentional islanding aims to split a network, by disconnecting lines,
into electrically-isolated islands. The challenge is that, if an
island is to be feasible, it must satisfy both static
constraints---load/generation balance, network constraints, system
limits---and dynamic constraints, \ie~for frequency and voltage
stability. Furthermore, the act of islanding must not cause a loss of
synchronism or voltage collapse. 

The majority of approaches to islanding aim to find, as a primary
objective, electromechanically stable islands. A popular approach
first uses slow coherency analysis to determine groupings of machines
with coherent oscillatory modes, and then aims to split the network
along the boundaries of these groups~\cite{YVY03,YVW04}.  Determining
the required cutset of lines involves considerations of
load-generation balance, power flows, and other constraints:
algorithms include pre-specification of boundaries~\cite{ASK+03},
exhaustive search~\cite{YVY03,YVW04}, minimal-flow minimal-cutset
determination using a combination of breadth- and depth-first
search~\cite{WV04}, graph simplification and partitioning~\cite{XV10},
spectral clustering~\cite{DGW+12}, and
meta-heuristics~\cite{LLC+09,AS12}. A key attraction of the
slow-coherence-based approach is that generator groupings are
dependent on machine properties and largely independent of fault
location and---to a lesser extent---operating point~\cite{YVW04}. If
the network can be split along the boundaries of these groups, while
not causing excessive load/generation imbalance or disruption, the
system is less likely to lose stability. Moreover, groupings and line
cutsets can be determined offline. As a consequence, the on-line
action of islanding is fast, and the approach has been demonstrated
effectively by simulations of real scenarios~\cite{YVH+07,XVM+11}.

Another approach uses ordered binary decision diagrams (OBDDs) to
determine balanced islands~\cite{SZL03}. Subsequently, power flow and
transient stability analyses can be used to iterate until feasible,
stable islands are found~\cite{SZL06}. In~\cite{JSS07},
a framework is proposed that iteratively identifies the controlling
group of machines and the contingencies that most severely impact
system stability. A heuristic method is used to search for a
splitting strategy that maintains a desired stability margin. Wang et
al.~\cite{WZH+10} employed a power flow tracing algorithm to first
determine the domain of each generator, \ie~the set of load buses that
``belong'' to each. Subsequently, the network is coarsely
split along domain intersections before refinement of boundaries to
minimize imbalances.

While it is known that the sensitivity of coherent machine groupings
to fault location is low, it is true that splitting the network along
the boundaries of \emph{a-priori} determined coherent groups is not,
in general, the only islanding solution that maintains
stability. Moreover, such islands may be undesirable in terms of other
criteria, such as the amount of load shed, the voltage profile or the
possibility that the impacted region may be contained within a larger
than necessary island. For example, in~\cite{DGW+12}, the
slow-coherence-based islanding of the $39$-bus New England system
isolates the network's largest generator in an island with no
load. In~\cite{TBG+13}, an optimization-based approach to islanding
and load shedding was proposed. A key feature is that, unlike many
other methods, it can take into account a 
part of the network that is
desired to be isolated---a troublesome area---when determining
islands, and isolate this while minimizing the expected amount of load
shed or lost. The problem is formulated as a single mixed integer
linear programming (MILP) problem, meaning that power balances, flows,
and operating limits may be handled explicitly when designing islands,
and satisfied in each island in a feasible solution.

The islanding MILP problem has similarities with the
\emph{transmission switching} problem~\cite{FOF08}, in that the
decision variables include which lines to disconnect, while
power flow constraints must be satisfied following any
disconnection. Both approaches---islanding and transmission
switching---may be seen as network topology optimization problems with
added power flow constraints. In both cases, inclusion of AC power
flow laws in the constraints results in a mixed integer nonlinear
program (MINLP), which is difficult to solve. Hence, linear DC power
flow has been used to date, resulting in a more computationally
favourable MILP or MIQP problem.

A disadvantage of the DC power flow model is that the effect of line
disconnections on network voltages is not considered. This is not
exclusive to MILP-based islanding and transmission switching; a number
of islanding approaches consider real power only, and assume that
reactive power may be compensated locally after
splitting. In~\cite{TBG+13}, however, cases were reported where a
solution could not be found to satisfy AC power flow and voltage
constraints when the islands were designed considering DC power flow,
\emph{even when sufficient reactive power generation capacity was
  present in each island}. Investigation found that local shortages or
surpluses of reactive power led to abnormal voltages in certain areas
of the network.

This paper presents a new method for controlled islanding that
respects voltage and reactive power constraints. A piecewise linear
approximation to AC power flow is developed and then used in a
MILP-based approach to islanding: decisions are which lines to
disconnect, which loads to shed and how to adjust generators. Results
on test networks show this eliminates the AC-infeasibilities reported
in~\cite{TBG+13}.  The method is flexible and able to deal with
different reasons for islanding.  For example, to minimize the load
shed while splitting the network so that coherent synchronous machines
remain in the same island.  Or, to split the network in two so as to
ensure that the most of it is left in a known safe state, isolated
from a troubled region that has been identified as a possible trigger
for cascading failures. The objective would be to minimize the load
that is planned to be shed, plus the expected extra load that might be
lost due to failures in the small island surrounding the troubled
region. There can be many reasons for suspecting trouble from a
region---\eg~incomplete or inconsistent measurements, estimates of
system stress such as closeness to instability or equipment operating
limits, indications of component failures, or other behaviour patterns
that simulations have shown to be correlated with cascading
failure~\cite{VBC+12}---but the precise definition of what evidence would lead to
islanding being initiated is complex and is beyond the scope of this
paper.

The organization of this paper is as follows. In the following
section, the piecewise linear AC power flow model is presented, and
its use is demonstrated in an Optimal Power Flow (OPF) problem. In
Section~\ref{sec:pwlislanding}, the islanding formulation is
described. Section~\ref{sec:results} presents computational results
for test networks. Conclusions are made in
Section~\ref{sec:conclusions}.

\section{Piecewise Linear AC Power Flow}
\label{sec:pwlacpower}

\subsection{A linear-plus-cosine model of AC power flow}
The linear ``DC'' model is a widely accepted approximation to AC power
flow, whose benefits (linearity, simplicity) often outweigh its
shortcomings. Recently, however, there has been renewed research
interest in the DC model itself~\cite{SJA09} and more accurate
alternative linearizations~\cite{CVB12}. Recent work~\cite{TBG+13} by
the authors found that a DC-based approach to controlled islanding
sometimes led to infeasible islands being created, mainly owing to
out-of-bound voltages and local shortages or surpluses of reactive
power. Motivated by this, this section presents a piecewise linear
approximation to AC power flow, in which voltage and reactive power
are modelled.

The AC power flow equations are described as follows.  Real and
reactive power balances at each bus $i \in \mathcal{B}$ give
\begin{align*}
  \sum_{g \in \mathcal{G}_i} p^{\ts{G}}_g &= \sum_{d \in \mathcal{D}_i} p^{\ts{D}}_d 
  + \sum_{l \in \mathcal{L}_i, j \in \mc{B}_l: j \neq i} p^{ij}_l 
  + G^{\ts{B}}_i v_i^2,\\
  \sum_{g \in \mathcal{G}_i} q^{\ts{G}}_g &= \sum_{d \in \mathcal{D}_i} q^{\ts{D}}_d 
  + \sum_{l \in \mathcal{L}_i, j \in \mc{B}_l: j \neq i} q^{ij}_l 
  - B^{\ts{B}}_i v_i^2,
\end{align*}
A line $l \in \mc{L}$ connects bus $i \in \mc{B}_l$ to bus
  $j\in\mc{B}_l, j \neq i$. The power flows from $i$ to $j$ are
\begin{align*}
  {p}^{ij}_l &=  v_{i}^2 G^{ii}_l + G^{ij}_l v_{i}v_{j} y_{ij} + B^{ij}_l v_{i}v_{j} z_{ij} ,\\
  {q}^{ij}_l &= -v_{i}^2 B^{ii}_l - B^{ij}_l v_{i}v_{j}y_{ij} + G^{ij}_l v_{i}v_{j} z_{ij},
\end{align*}
with a similar expression from $j$ to $i$, where
\begin{gather*}
\tau_l^2 G_{l}^{ii} = G_{l}^{jj} = -\tau_lG_l^{ij} = -\tau_lG_l^{ji} = g_l,\\
\tau_l^2 B_{l}^{ii} = B_{l}^{jj} = -\tau_lB_l^{ij} = -\tau_lB_l^{ji} = b_l + 0.5b^\ts{C}_l.
\end{gather*}
The convention is for a transformer to be located at the \emph{from}
end (bus $i$) of a branch.

The standard ``DC'' approximation to AC power flow linearizes these
equations by using the approximations $v_{i} = v_{j} = 1$, $z_{ij}
=\theta_{ij}$, $y_{ij} = 1$, and $b_l \gg g_l \approx 0$ yielding
${p}^{ij}_l = B^{ij}_l \theta_{ij}$. The reactive power variables and
equations are dropped. In the model in this paper, voltage and
reactive power are retained. Expanding the line flows about $v_{i} =
1$, $v_{j} = 1$ and $\theta_{ij} = 0$ (hence $y_{ij} = 1, z_{ij} =
0$):
\begin{align*}
p^{ij}_l &\approx G^{ii}_l (2v_{i} - 1) + G^{ij}_l\bigl( v_{i} + v_{j} + y_{ij} - 2 \bigr) + B^{ij}_l z_{ij},\\
q^{ij}_l &\approx B^{ii}_l (1 - 2v_{i}) - B^{ij}_l\bigl( v_{i} + v_{j} + y_{ij} - 2 \bigr) + G^{ij}_l z_{ij}.
\end{align*}
In a standard linearization, the small-angle approximations would then
be used: $y_{ij} = \cos\theta_{ij} \approx 1$ and $z_{ij} = \sin \theta_{ij}
\approx \theta_{ij}$. %
%
%
Tab.~\ref{tab:approxerrors} gives the maximum absolute errors for each
of the constituent terms in the linearized flows, over a typical range
of operating voltages and angles, \ie~$0.95 \leq v_{i} \leq 1.05$ at
each end of the line, and $|\theta_{ij}| \leq \unit{40}{\degree}$. The
cosine approximation incurs the largest error. Fig.~\ref{fig:PQerrors}
shows maximum and minimum power flows and errors over this range of
voltages and angles for a line with $g_l = 1, b_l = -5, b^\ts{C}_l =
1$. Approximation errors are obtained for when the $y_{ij} =
\cos\theta_{ij}$ term is approximated as $1$ (a linear model) and
modelled exactly (linear plus cosine). In both cases, $z_{ij} =
\theta_{ij} \approx \sin \theta_{ij}$ is used. Although little
reduction in errors is apparent in the real flows, the importance of
modelling the cosine term is clear for reactive flows.
\begin{table}[b]
\centering\footnotesize
\caption{Approximation errors in line flow terms}
\begin{tabular}{ccr}
\toprule
Term & Approximation & Max abs error \\
\midrule
$v_{i}^2$ & $2v_{i} - 1$ & $0.0025$ \\
$v_{i} v_{j} y_{ij}$ & $v_{i} + v_{i} + y_{ij} - 2$ & $0.0253$ \\
$v_{i} v_{j} z_{ij}$ & $z_{ij}$ & $0.0659$ \\
$y_{ij}$ & $1$ & $0.2340$ \\
$z_{ij}$ & $\theta_{ij}$ & $0.0553$ \\
\bottomrule
\end{tabular}
\label{tab:approxerrors}
\end{table}
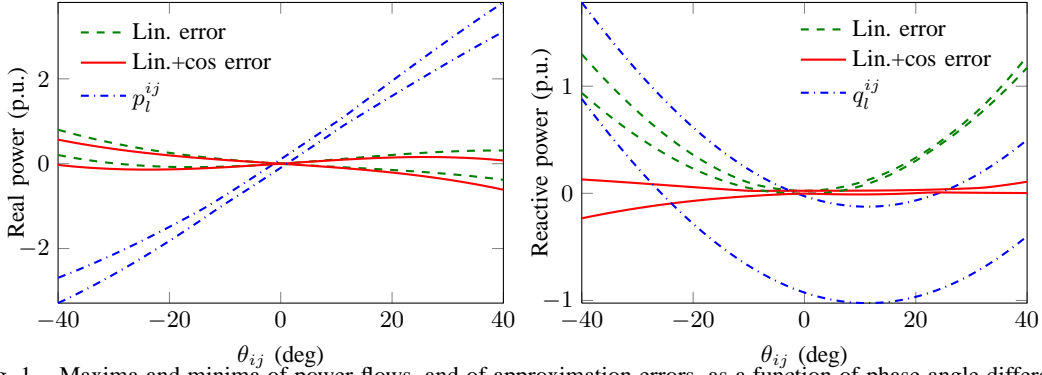
\begin{figure}[t]
\centering\footnotesize
\input{trodd1.tikz}
\vspace{-2em}
\caption{Maxima and minima of power flows, and of approximation errors, as a
  function of phase angle difference.}
\label{fig:PQerrors}
\end{figure}

A similar analysis shows that including the sine term (instead of its
linearization) in addition to the cosine term reduces the error in the
real flows slightly, but makes no significant difference to the
reactive power. Since the infeasibilities that occur using the DC
approach to islanding are mainly owing to the reactive power and
voltage limits~\cite{TBG+13}, the appropriate approximation to use is
the linear-plus-cosine one. And although cosine terms cannot be used
directly in an MILP model, they can be modelled to arbitrary levels of
accuracy by piecewise linear functions. The next section demonstrates
the use of the model in an OPF formulation.

\subsection{Piecewise linear AC OPF}

The piecewise linear (PWL) AC OPF problem is defined as
\begin{equation*}
\min \sum_{g \in \mc{G}} c_g \bigl(p^\ts{G}_g\bigr)
\end{equation*}
subject to, $\forall i \in \mc{B}$, the linearized power balances:
\begin{subequations}
  \begin{align}
	  \sum_{g \in \mc{G}_i} p^{\ts{G}}_g &= \sum_{d \in \mathcal{D}_i} P^{\ts{D}}_d 
	  + \sum_{l \in \mathcal{L}_i, j \in \mc{B}_l: j \neq i} p^{ij}_l 
	  + G^{\ts{B}}_i(2v_i-1),\label{eq:linKCL1}\\
	  \sum_{g \in \mathcal{G}_i} q^{\ts{G}}_g &= \sum_{d \in \mathcal{D}_i} Q^{\ts{D}}_d 
	  + \sum_{l \in \mathcal{L}_i, j \in \mc{B}_l: j \neq i} q^{ij}_l 
	  - B^{\ts{B}}_i (2v_i-1),\label{eq:linKCL2}
\end{align}
\end{subequations}
Line flows for all $l \in \mc{L}, i, j \in \mc{B}_l: i \neq j$:
	\begin{align*}
		{p}^{ij}_l &= G^{ii}_l (2v_{i} - 1) + G^{ij}_l\bigl( v_{i} + v_{j} + y_{ij} - 2 \bigr) + B^{ij}_l \theta_{ij},\\
		{q}^{ij}_l &= B^{ii}_l (1 - 2v_{i}) - B^{ij}_l\bigl( v_{i} + v_{j} + y_{ij} - 2 \bigr) + G^{ij}_l \theta_{ij}.
	\end{align*}
The $N$-piece PWL approximation to $\cos\theta_{ij}$. For all $l \in \mc{L}, i, j \in \mc{B}_l: i \neq j$.
  \begin{equation}
  y_{ij} =  h_{ij,k} \theta_{ij} + d_{ij,k}, \forall \theta_{ij} \in [x_{ij,k} , x_{ij,k+1}], k=0\ldots N-1,
  \label{eq:xidef1}%
\end{equation}
where $h_{ij,k}$ and $d_{ij,k}$ are chosen so that the
  approximation coincides with $\cos x$ at breakpoints
$\{x_{ij,0},\ldots,x_{ij,N}\}$.
System limits are applied:
\begin{gather}
  V^-_i \leq v_i \leq V^+_i, \forall i \in \mc{B},\notag\\
  P^{{\ts{G}-}}_g \leq p^\ts{G}_g \leq P^{{\ts{G}+}}_g, \forall g \in \mc{G},\notag\\
  Q^{{\ts{G}-}}_g \leq q^\ts{G}_g \leq Q^{{\ts{G}+}}_g, \forall g \in \mc{G},\notag\\
	p_l^{ij} + p_l^{ji} \leq P_l^{\ts{L}+}, \forall l \in \mc{L}.\label{eq:linelim}
\end{gather}
Note that line flow limits are limits on real power ($I^2R$) loss. If
an MVA limit $S^{\ts{L}+}_l$ is given, this may be converted by
assuming nominal voltage, \ie~$P_l^{\ts{L}+} = \frac{g_l}{g^2_l +
  b^2_l} \bigl(S^{\ts{L}+}_l\bigr)^2$.

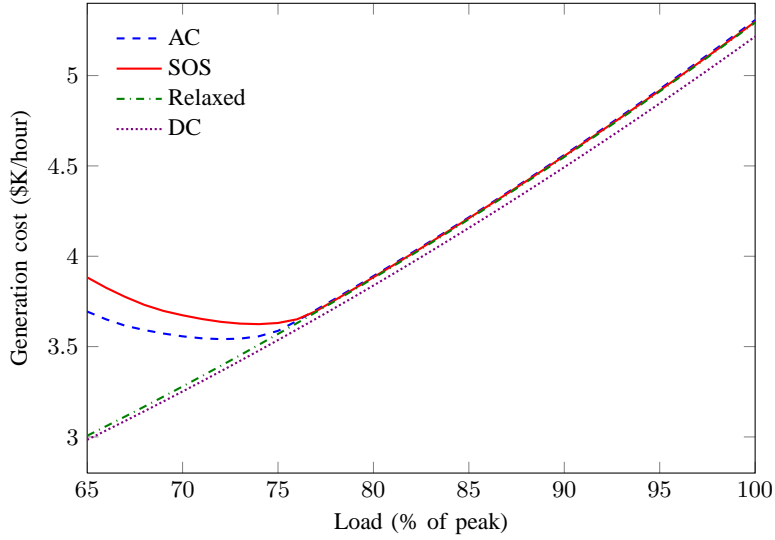
\begin{figure}[t]
\centering\footnotesize
\input{trodd2KM.tikz}
\vspace{-1em}
\caption{Generation costs as a function of load for the $9$-bus network.}
\label{fig:9bus}
\end{figure}

The implementation of the PWL model of
$\cos\theta_{ij}$~\eqref{eq:xidef1} requires either binary variables
or special ordered sets of type~2 (SOS-2)~\cite{BF76}. The overall
problem is then, depending on $c_g$, a mixed integer linear or
quadratic program (MILP or MIQP). If ~\eqref{eq:xidef1} is replaced by
its relaxation $y_{ij} \leq h_{ij,k} \theta_{ij} + d_{ij,k}$, then the
problem becomes a convex optimization problem and no binary variables
or SOS sets are needed. Since real and reactive line losses decrease
as $y_{ij}$ increases, it is tempting to assume that equality will
hold for one of the PWL sections, and this relaxation will yield a
tight result. However, as Fig.~\ref{fig:9bus} shows, situations exist
where the SOS formulation is necessary. This shows optimal generation
costs against load level, as obtained by OPFs using AC, PWL with SOS,
relaxed PWL, and DC power flow models. The network is the WSCC $9$-bus
network modified to set voltage limits to $\pm 5\%$ and the lower
reactive power limit for each generator is raised from $-300$ to $-5$
Mvar. This means that at low load levels the generators find it
increasingly difficult to balance the reactive power, as more lines
become sources rather than sinks of reactive power, and the generation
cost rises with falling load.
While the SOS PWL is able to capture this effect, the relaxed PWL and
DC-based models are not; the former ``cheats'' by having some lines
continue to store reactive power irrespective of their end voltages
and angles---allowed because $y_{ij} < h_{ij,k}\theta_{ij} + d_{ij,k}$
is permitted---and this allows more of the real power to be generated
by the cheaper generators.

\section{A formulation for system islanding using piecewise linear AC power flow}
\label{sec:pwlislanding}

In~\cite{TBG+13}, the problem of determining how to split a
transmission network into islands is considered. The aim is to limit
the effects of possible cascading failures and prevent the onset of
wide-area blackouts by re-configuring the network---\emph{via} line
switching---so that problem areas are isolated. The MILP-based method
defines two \emph{sections} of the network. All of the buses that must
be isolated are pre-assigned to section $0$, and the optimization
determines which other buses and lines to place in section~$0$. All
the remaining components are in section~$1$. This creates at least two
islands. The optimization will also determine the best strategy to
adjust generation and shed load so as to establish a load-generation
balance in each island while respecting all network equations and
operating constraints after the split.

\subsection{Motivation: effect of topology changes on voltage profile}
\label{sec:infeas_example}

Solution of the MILP islanding problem provides a set of lines to
switch, loads to shed and generators to adjust. However, if only
the DC power flow equations are included in the constraints, the
effects of changing the network topology on voltages and reactive
power flows is not considered. Thus, in~\cite{TBG+13}, an AC optimal
load shedding (OLS) problem is solved \emph{after} the MILP islanding
problem, using the islanded network topology. If a solution to this
can be found, the islanded network is feasible with respect to AC
power flow and operating constraints. The solution provides the
correct generator output and load adjustments to make, now having
considered voltage and reactive power.

However a number of the islanding solutions in~\cite{TBG+13} were AC
infeasible, primarily due to violation of voltage bounds;
solutions could be recovered
by relaxing the normal limits.
\begin{figure}[t]
\centering\scriptsize
\vspace{-1em}
\subfloat{%
\psfrag{1}[t][t]{$1$}%
\psfrag{2}[t][t]{$2$}%
\psfrag{3}[t][t]{$3$}%
\psfrag{4}[t][t]{$4$}%
\psfrag{5}[t][t]{$5$}%
\psfrag{6}[t][t]{$6$}%
\psfrag{7}[t][t]{$7$}%
\psfrag{8}[t][t]{$8$}%
\psfrag{9}[t][t]{$9$}%
\psfrag{10}[t][t]{$10$}%
\psfrag{11}[t][t]{$11$}%
\psfrag{12}[t][t]{$12$}%
\psfrag{13}[t][t]{$13$}%
\psfrag{14}[t][t]{$14$}%
\psfrag{15}[t][t]{$15$}%
\psfrag{16}[t][t]{$16$}%
\psfrag{17}[t][t]{$17$}%
\psfrag{18}[t][t]{$18$}%
\psfrag{19}[t][t]{$19$}%
\psfrag{20}[t][t]{$20$}%
\psfrag{21}[t][t]{$21$}%
\psfrag{22}[t][t]{$22$}%
\psfrag{23}[t][t]{$23$}%
\psfrag{24}[t][t]{$24$}%
\psfrag{x}[t][t]{$\sim$}%
\includegraphics[width=0.315\textwidth,height=7.5cm]{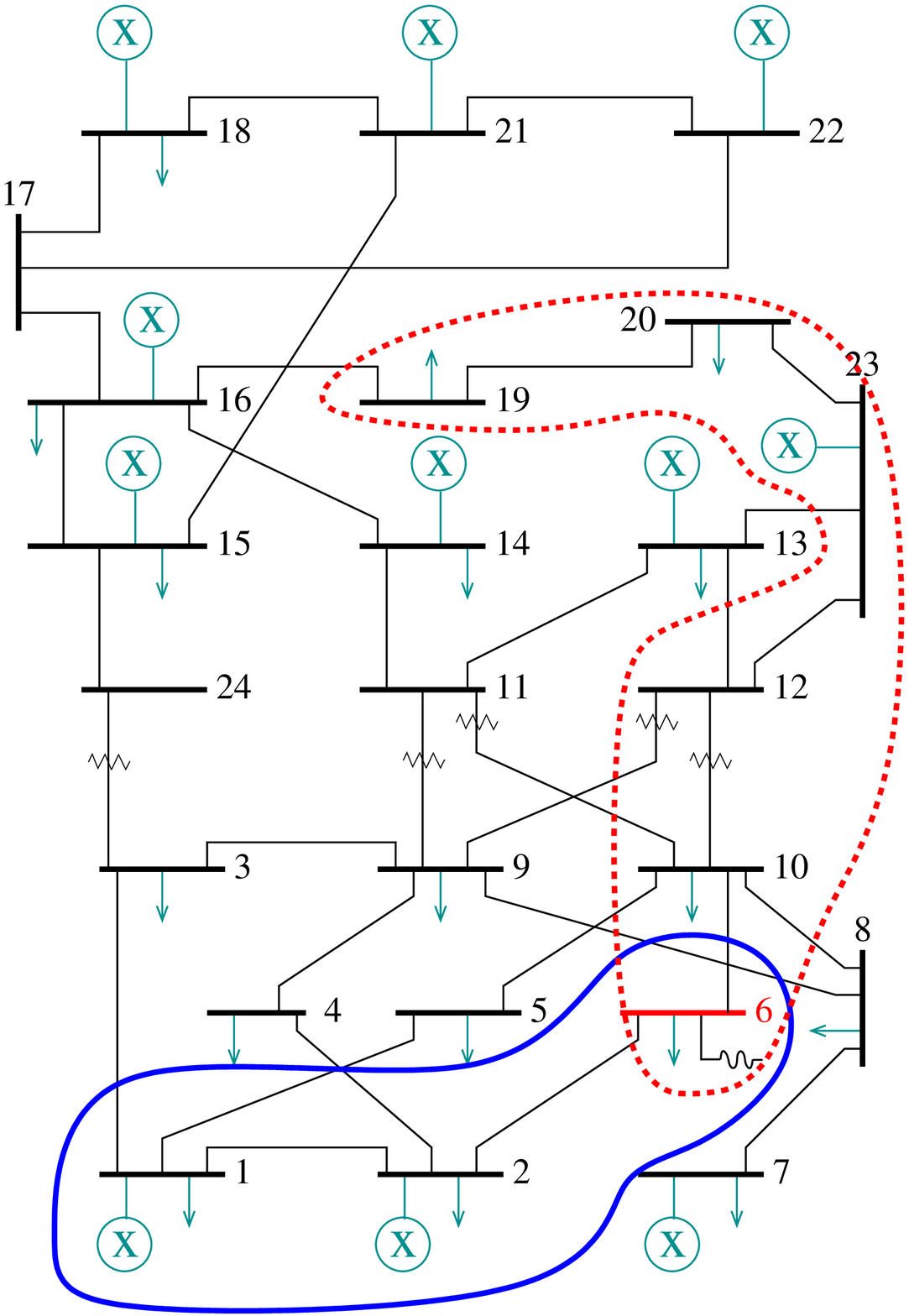}
}\hfil
\subfloat{%
\psfrag{1}[t][t]{$1$}%
\psfrag{2}[t][t]{$2$}%
\psfrag{3}[t][t]{$3$}%
\psfrag{4}[t][t]{$4$}%
\psfrag{5}[t][t]{$5$}%
\psfrag{6}[t][t]{$6$}%
\psfrag{7}[t][t]{$7$}%
\psfrag{8}[t][t]{$8$}%
\psfrag{9}[t][t]{$9$}%
\psfrag{10}[t][t]{$10$}%
\psfrag{11}[t][t]{$11$}%
\psfrag{12}[t][t]{$12$}%
\psfrag{13}[t][t]{$13$}%
\psfrag{14}[t][t]{$14$}%
\psfrag{15}[t][t]{$15$}%
\psfrag{16}[t][t]{$16$}%
\psfrag{17}[t][t]{$17$}%
\psfrag{18}[t][t]{$18$}%
\psfrag{19}[t][t]{$19$}%
\psfrag{20}[t][t]{$20$}%
\psfrag{21}[t][t]{$21$}%
\psfrag{22}[t][t]{$22$}%
\psfrag{23}[t][t]{$23$}%
\psfrag{24}[t][t]{$24$}%
\psfrag{x}[t][t]{$\sim$}%
\includegraphics[width=0.315\textwidth,height=7.5cm]{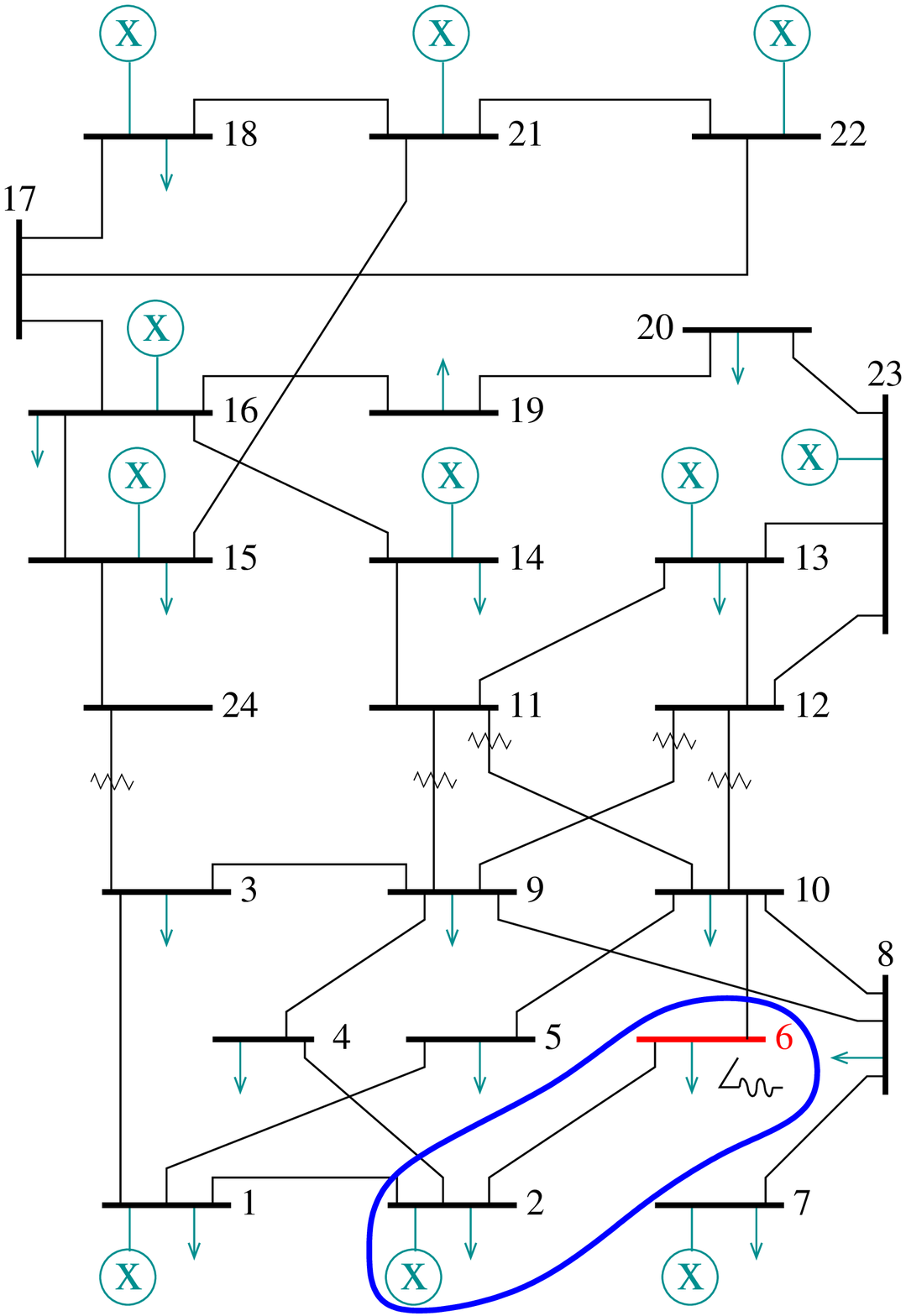}
}
\caption{IEEE 24-bus RTS with bus $6$ isolated. On the left, DC method
  (solid) and PWL method without shunt switching (dashed). On the
  right, PWL method with shunt switching.}
\label{fig:24bus}%
\end{figure}
One such example, for the $24$-bus IEEE Reliability Test System
(RTS)~\cite{IEEERTS96}, is described as follows. Given the problem of
isolating bus $6$ while minimizing the expected load shed or lost, the
optimal solution islands buses $1$, $2$ and $6$, as indicated in
Fig.~\ref{fig:24bus}. There remains sufficient real power capacity in
both islands to meet demand, and no load is shed. Moreover, but not by
design, there is sufficient reactive power capacity in each island to
meet the total reactive power demand. Despite this, a feasible
solution to the AC-OLS cannot be found. Softening the voltage
bounds recovers a solution, but with an abnormally low voltage of
$0.6443$~p.u. at bus $6$ and an over-limit flow on line
$(2,6)$. Further inspection reveals that this situation has arisen
because of the disconnection of line $(6,10)$, a cable with high shunt
capacitance. The passive shunt reactor at bus $6$ would, in normal
circumstances, balance locally the excess reactive power and maintain
a satisfactory voltage profile. This problem could be avoided by
linking together the disconnection of line $(6,10)$ and the shunt
reactor at $6$. The optimal solution when these actions are linked is
shown in the right-hand diagram of Fig.~\ref{fig:24bus}, and it yields
a better feasible solution than when the reactor is not
disconnected. Rules like this are easy to incorporate in the model;
however, it is difficult \emph{a-priori} to define all possible
rules. A better approach is to allow the model to decide the
combination of equipment to disconnect, and when this is done the
optimal solution disconnects both line $(6,10)$ and the reactor
at bus $6$, giving the right-hand Fig.~\ref{fig:24bus} solution.
The models and result for this example are given in
Sections~\ref{OtherCompSwitching} and ~\ref{CompSwitchResults}.

This is just one example of where an islanding solution formed by
considering only real power---even if network constraints are
included---is unsatisfactory. It also shows that even if reactive
power balance is achieved within each island, local shortages or
surpluses can lead to an abnormal voltage profile. Many test networks
are prone to this problem~\cite{TBG+13}. Moreover, it
is not just system islanding that is susceptible; DC-based
transmission switching also does not consider the consequences on
voltage of disconnecting lines. Thus, there is a need for network
topology optimization methods that can determine AC-feasible
solutions, but without having to resort to solving the full MINLP
problem. The focus of this paper is topology optimization for the
purpose of islanding, and in the next section, a formulation is
presented that uses the PWL model of AC power flow.
 
\subsection{Formulation of constraints for islanding}

The problem is to decide which lines to switch in order to isolate a
part of the network. Separation of sections is enforced by sectioning
constraints. The islanded network must satisfy power balance and flow
equations and operating limits, and so these are
included as constraints in the problem.

\subsubsection{Sectioning constraints}

Define $\mc{B}^0$ and $\mc{B}^1$, where $\mc{B}^0 \cap \mc{B}^1 =
\emptyset$, as the subsets of buses that are desired to be
separated. For now, the motivation for this separation is left open,
but it may be that these buses in, say, $\mc{B}^0$ represent a failing
area of the network, or are associated with a coherent group of
synchronous machines that will be separated from other groups. The
proposed approach will split the network into two \emph{sections}:
section~$0$ will contain all buses in $\mc{B}^0$ and section~$1$ all
buses in $\mc{B}^1$. For a bus $i \in \mc{B}$, $\gamma_i$ denotes the
section ($0$ or $1$) to which that bus is assigned. That is, if $i$ is
to be placed in section $0$, then $\gamma_i = 0$. Separation between
sections is achieved by switching lines: $\rho_l$ denotes the
connection status of a line $l$, and the convention followed is for
$\rho_l=0$ when $l$ is disconnected. The exact boundaries of each
section will depend on the objective, defined later, and the
optimization will determine how to assign to sections those buses not
in $\mc{B}^0$ or $\mc{B}^1$, in order achieve balance and optimize the
objective. However, the following constraints enforce the separation
of sections $0$ and $1$, without defining precisely their boundaries.
\begin{subequations}
\begin{align}
\rho_l &\leq 1 + \gamma_{i} - \gamma_{j}, 
\forall l \in \mc{L}, i, j \in \mc{B}_l: i\neq j,\label{eq:sec1a}\\
\gamma_i &= s, \forall i \in \mc{B}^s, s \in \{0,1\}.\label{eq:sec3a}
\end{align}
\label{eq:IPsectioning}%
\end{subequations}
\subsubsection{Power flow}

The remainder of the constraints are concerned with achieving a
balanced, steady state for the islanded network. It is assumed that
generators are permitted to make only small-scale changes to output or
be switched off, and loads may be fully or partly shed in order to
maintain a balance. As a consequence of these changes and the
topological changes, bus voltages, angles and line flows will change,
and so must be modelled to ensure satisfaction of network constraints
and operating limits.

First, the power balances, \eqref{eq:linKCL1} and~\eqref{eq:linKCL2},
are included without modification. Next, the line flow equations are
modified so that when a line is disconnected, power flows across it
are zero irrespective of its end bus voltages and angles. To assist
this, we introduce \emph{line} variables---$v_l^i$ and $v_l^j$ as end
voltages and $\theta^{ij}_l$ as the angle difference---that are
distinct from \emph{bus} variables $v_i$, $v_j$ and $\theta_{ij}$. The
following constraints control the relationship between line variables
and bus variables. For a line $l \in \mc{L}$ with end buses $i$ and
$j$,
\begin{subequations}
\begin{gather}
-\Theta_l \rho_l \leq \theta^{ij}_l \leq \Theta_l \rho_l,\label{eq:thetaswitch1}\\
-\Theta_l^+(1-\rho_l) \leq \theta^{ij}_l - \theta_{ij} \leq \Theta_l^+(1-\rho_l),\label{eq:thetaswitch2}
\intertext{$\forall i \in \mc{B}_l:$}
0 \leq v_{i} - v^i_l \leq ( V^+_{i} - V^-_{i} ) (1-\rho_l), \label{eq:vswitch1}\\
V^-_{i} \leq v^i_l \leq V^-_{i} + (V^+_{i} - V^-_{i})\rho_l,\label{eq:vswitch3}
\intertext{and $\forall i \in \mc{B}$,}
V^-_i \leq v_i \leq V^+_i,\label{eq:vbounds}
\end{gather}
\label{eq:flow1}
\end{subequations}
where $\Theta_l^+ \geq \Theta_l$~is a ``big-$M$'' constant. Of these, \eqref{eq:thetaswitch1}
and~\eqref{eq:thetaswitch2} force equality of $\theta^{ij}_l$ and
$\theta_{ij}=\delta_{i} - \delta_{ij}$ for a connected line, but set $\theta^{ij}_l
= 0$ for a disconnected line while allowing the \emph{bus} angles
$\delta_{i}$ and $\delta_{j}$ to vary independently. Likewise, if
$\rho_l = 1$ then, by~\eqref{eq:vswitch1}, $v^i_l = v_{i}$ and $v^j_l
= v_{j}$. However, if $\rho_l = 0$ then the line voltages are set to
minimum values---$v^i_l = V^-_{i}$ and $v^j_l = V^-_{j}$---independent
of the bus voltages $v_i$ and $v_j$. 

This switching between line and bus variables is made use of in
modified line flow equations. For a line $l$,
\begin{subequations}
\begin{align}
\begin{split}
	p^{ij}_l & = G^{ii}_l (2v^{i}_{l} - 1) + G^{ij}_l\bigl( v^{i}_{l} + v^{j}_{l} + y^{ij}_l - 2 \bigr) + B^{ij}_l \theta^{ij}_l\\
&\quad - \bigl(G^{ii}_l(2 V^-_{i} - 1) +
G^{ij}_l(V^-_{i}+V^-_{j}-1)\bigr)(1-\rho_l),
\end{split}\raisetag{2.3em}\\
\begin{split}
	q^{ij}_l & = B^{ii}_l (1 - 2v^{i}_{l}) - B^{ij}_l\bigl( v^{i}_{l} + v^{j}_{l} + y^{ij}_l - 2 \bigr) + G^{ij}_l \theta^{ij}_l \\ 
&\quad-\bigl( B^{ii}_l(1-2 V^-_{i}) -B^{ij}_l(V^-_{i}+V^-_{j}-1) \bigr) (1-\rho_l),
\end{split}
\raisetag{2.3em}
\end{align}
\label{eq:flow2}
\end{subequations}
and $y^{ij}_l$ is given by~\eqref{eq:xidef1}, using
$\theta^{ij}_l$. Note that since $\theta^{ij}_l = 0$ if $\rho_l=0$,
then $y^{ij}_l=1$ for a disconnected line. Hence, if $\rho_l = 0$ then
$p^{ij}_l=0$, irrespective of $v_i$, $v_j$ and $\theta_{ij}=\delta_i -
\delta_j$. If $\rho_l = 1$, the normal power flow equations are
recovered.

\subsubsection{Operating constraints}
\label{sec:opercons}

In the short time available when islanding in response to a
contingency, any extra generation that is needed will be achieved by a
combination of the ramping-up of on-line units and the commitment of
fast-start units. For simplicity, fast-start units are not considered
in the examples in this paper. We assume that a generator that is
operating can either have its input mechanical power disconnected, in
which case real output power drops to zero in steady state, or its
output can be set to a new value within a small interval,
$\bigl[P^{\ts{G}-}_g, P^{\ts{G}+}_g\bigr]$, say, for generator $g$,
around the pre-islanded value. The limits will depend on the ramp and
output limits of the generator, and the amount of immediate or
short-term reserve capacity available to the generator. For the test
scenarios in Section~\ref{sec:results}, a time limit of $2$ minutes is
assumed for ramping, but the formulation permits any choice. This
choice should be informed by existing post-contingency response
protocols.  For reactive power, it is assumed that a new output can be
set in some range $Q^{\ts{G}-}_g$ to $Q^{\ts{G}+}_g$. The set of
possible real and reactive power outputs of a generator is usually
convex.  For the test scenarios in Section~\ref{sec:results}, the
bounds on the real and reactive power are independent.  In the more
general case, since the range of values for the real power output is
small, the feasible region for the problem is a narrow slice through a
convex set, and---except when the real power output is close to its
upper limit---it is a good approximation to treat the real and
reactive power bounds as independent. If this is not the case, it is
straightforward to add constraints that couple $p^\ts{G}_g$ and
$q^\ts{G}_g$.

The operating regime is modelled by the constraints
\begin{subequations}
\begin{gather}
\zeta_g P^{\ts{G}-}_g \leq p^{\ts{G}}_g \leq \zeta_g P^{\ts{G}+}_g, \forall g \in \mc{G},\\
Q^{\ts{G}-}_g \leq q^{\ts{G}}_g \leq Q^{\ts{G}+}_g, \forall g \in \mc{G},\\
  \zeta_g = 1, \forall g \in \bigl\{\mc{G}: P^{\ts{G}-}_g = 0 \bigr\} \cup \mc{G}^1.
\end{gather}%
\label{eq:gencons}%
\end{subequations}%
Here, $\zeta_g$ is a binary variable and denotes the on/off setting of
the real power output, and $\mc{G}^1$ is a subset of generators which
are required to remain on.

For loads, because of the limits on generator outputs and network
constraints, it may not be possible after islanding to fully supply
all loads. It is therefore assumed that some shedding of loads is
permissible. Note that this is \emph{intentional} shedding, not
automatic shedding as a result of low voltages or frequency. To
implement this in the real network there has to be central control
over equipment. For all $d \in \mc{D}$,
\begin{subequations}
\begin{align}
  p^{\ts{D}}_d &= \alpha_d P^{\ts{D}}_d,\label{eq:loadshedP}\\
  q^{\ts{D}}_d &= \alpha_d Q^{\ts{D}}_d,\label{eq:loadshedQ}
\end{align}
\label{eq:loadcons}
\end{subequations}
where $0 \leq \alpha_d \leq 1$.

Finally, line limits are applied via constraint~\eqref{eq:linelim}.

\subsubsection{Rules for other component switching} 
\label{OtherCompSwitching}

As motivated by Section~\ref{sec:infeas_example}, sometimes it is
necessary to have rules for switching components or adjusting controls
in different situations. Such rules can easily be included in the
formulation using standard techniques for deriving constraints from
logical rules~\cite{Williamsbook}. For example, the switching of a
shunt component at a bus $i$ can be modelled by introducing binary and
continuous variables, $\xi_i$ and $u_i$ respectively, constraints
\begin{gather*}
\xi_i (2V^-_i -1) \leq u_i \leq \xi_i (2V^+_i -1 ),\\
-(1-\xi_i) (2V^-_i -1) \leq u_i - (2v_i -1) \leq (1-\xi_i) (2V^+_i -1 ),
\end{gather*}
and replacing the $G^\ts{B}_i(2v_i-1)$, $B^\ts{B}_i(2v_i-1)$ terms
in~\eqref{eq:linKCL1} and \eqref{eq:linKCL2}
with $G^\ts{B}_i u_i$ and $B^\ts{B}_i
u_i$, respectively. In Section~\ref{CompSwitchResults} this is explored
further for the $24$-bus example.

\subsection{Objective functions for islanding}

The general aim is to split the network, separating the two sections
$0$ and $1$, yet leaving it in a feasible state of operation. The
specific motivations and objectives for islanding are discussed in
this section. Clearly, if a network can be partitioned with minimal
disruption to load, and with minimal disturbances to generators, then
its chances of viable operation until future restoration are
increased.

\subsubsection{Isolating uncertain regions and maximizing expected load supply}

We assume that there is an identifiable localized area of the network
that is believed could be a trigger for cascading failure. Similar to
the approach in~\cite{TBG+13}, the goal is to include this area of
potential trouble in an island, leaving the rest of the network in a
known, secure steady state.  The sets $\mc{B}^0$ and $\mc{L}^0$
consist of all buses and lines in the troubled area and, additionally,
any buses and lines whose status is uncertain.  To ensure section $1$
contains no uncertain components, all lines $l \in \mc{L}^\ts{0}$ remaining
in this section are disconnected by replacing ~\eqref{eq:sec1a} by
\begin{equation}
\rho_l \leq 1 - \gamma_{i}, \forall i \in \mc{B}_l.\label{eq:sec2a}
\end{equation}

Because section $0$ may contain failing components or be in an
uncertain state, it is assumed there is a risk of not being able to
supply any load placed in that section.
Accordingly, a load loss penalty
$0\leq\beta_d<1$ is defined for a load $d$, which may be interpreted
as the probability of being able to supply a load if placed in section
$0$. Suppose a reward $R_d$ is obtained per unit supply of load
$d$. If $d$ is placed in section $1$ a reward $R_d$ is realized per
unit supply; however, if $d$ is placed in section $0$, a lower reward
of $\beta_dR_d<R_d$ is realized.

The objective is then to
maximize the expected total value of load supplied:
\begin{equation}
J^\ts{exp load} = \sum_{d \in \mc{D}} R_d P^\ts{D}_d (\beta_d \alpha_{0d} + \alpha_{1d}),\label{eq:obj1}
\end{equation}
where $\alpha_d = \alpha_{0d} + \alpha_{1d}$, and $0 \leq \alpha_{1d}
\leq \gamma_b, \forall b \in \mc{B}, d \in \mathcal{D}_b$. Here a new
variable $\alpha_{sd}$ is introduced for the load $d$ delivered in
section $s \in \{0,1\}$. If $\gamma_b = 0$ (and so the load at bus $b$
is in section~$0$), then $\alpha_{1d} = 0$, $\alpha_{0d} = \alpha_d$,
and the reward is $\beta_d R_d P^\ts{D}_d \alpha_d$. Conversely, if
$\gamma_b=1$ then $\alpha_{1d}=\alpha_d$ and $\alpha_{0d}=0$, giving a
larger reward $R_dP_d^\ts{D}\alpha_d$. Thus,
maximizing~\eqref{eq:obj1} gives a preference for $\gamma_b=1$ and a
smaller section $0$, so that the impacted area is limited.

\subsubsection{Promoting generator coherency}

Another aim is to ensure the synchronicity of generators within
islands. Large disturbances in the network cause electro-mechanical
oscillations, which can lead to a loss of synchronism. A popular
approach is to split the system along boundaries of near-coherent
generator groups, as determined by slow-coherency
analysis~\cite{AKW+80}. Thus, weak connections between
machines---which give rise to slow, lightly-damped oscillations---are
cut, leaving separate networks of tightly-coupled, coherent machines.

Consider those buses in the network with generators attached, the set
of which is defined as $\mathcal{B}^\ts{G}$, and define
$\mathcal{B}^\ts{GG} \triangleq \bigl\{ (i,j)\in
\mc{B}^\ts{G}\times\mc{B}^\ts{G} : j > i \bigr\}$ as the set of all
pairs of such buses. For what follows, it may be assumed that multiple
units at a bus are tightly coupled and are aggregated to a single
unit. The dynamic coupling, $W_{ij}$, between a pair of machines at
buses $(i,j) \in \mathcal{B}^\ts{GG}$ may be determined from
slow-coherency analysis. For example, assuming as in~\cite{DGW+12} the
undamped second order swing equation,
\begin{equation*}
W_{ij} = \pd{(\dot{\omega}_i - \dot{\omega}_{j})}{(\delta_{i}-\delta_{j})} = \biggl( \frac{1}{M_i} + \frac{1}{M_{j}} \biggr) \pd{P_{ij}}{\delta_{ij}},
\end{equation*}
where $M_i$, $\omega_i$, $\delta_i$ are the inertia constant, angular
frequency and rotor angle of the machine at bus $i$, and
$\pd{P_{ij}}{\delta_{ij}}$ is the synchronizing power coefficient or
``stiffness'' between machines at $i$ and $j$. To favour, in the
objective, separating loosely-coupled generators, introduce a new
variable $0\leq \eta_{ij} \leq 1$ for all $(i,j) \in
\mc{B}^\ts{GG}$. Then the constraint
\begin{equation}
-\eta_{ij} \leq \gamma_i - \gamma_j \leq \eta_{ij},
\end{equation}
sets $\eta_{ij}$ to $1$ if generator buses $i$ and $j$ are in
different sections of the network (and hence electrically isolated),
but otherwise may be zero. Minimizing the function
\begin{equation}
J^\ts{coh} = \sum_{(i,j) \in \mc{B}^\ts{GG}} W_{ij} \eta_{ij} \label{eq:obj2}
\end{equation}
gives a preference for machines in different sections having small
$W_{ij}$, \ie~being weakly coupled, and those within the same section
have stronger coupling. This may be used in conjunction
with~\eqref{eq:obj1}, \ie~$\max J^\ts{exp load} - kJ^\ts{coh}$,
with weighting $k > 0$, so that section $0$ is the ``unhealthy''
section, and the expected load supply is maximized while keeping
together strongly-coupled machines.

Minimizing~\eqref{eq:obj2} alone will favour keeping all machines in
the same section, and to force the machines apart additional
constraints may be needed. Alternatively, the following implementation
splits the network directly into coherent groups, making different use
of the sets $\mc{B}^0$ and $\mc{B}^1$.

\subsubsection{Splitting into coherent groups}

Suppose that coherent groups of generators have been determined, and
that assigned to $\mc{B}^0$ and $\mc{B}^\ts{1}$ are those buses in
$\mc{B}^\ts{G}$ corresponding to machines in different groups. For
example, $\mc{B}^0$ may contain the critical coherent group of machines, and
$\mc{B}^1$ all others. The sectioning constraints will ensure that the
machines are separated, but which other buses are assigned to each
section is determined by the optimization. The solution that minimizes
the amount of load shed can be found by maximizing the function
\begin{equation}
J^\ts{load} = \sum_{d \in \mc{D}} \alpha_d P^\ts{D}_d.
\label{eq:sumobj}
\end{equation}
Alternatively, to seek a solution that changes the generator outputs
the minimally from their initial values $P^{\ts{G}0}_g$, minimize
\begin{equation}
J^\ts{gen} = \sum_{g \in \mc{G}} t_g
\label{eq:genmoveobj}%
\end{equation}
where $t_g \geq 0$, $t_g \geq p^\ts{G}_g - P^{\ts{G}0}_g$, and $t_g
\geq -p^\ts{G}_g + P^{\ts{G}0}_g, \forall g \in \mc{G}$. The
sectioning constraints ensure that the machines are split into two
sections. If further separation is required, the optimization can be
re-run on each island of the network.

\subsubsection{Penalties} 
Often there may be multiple feasible solutions with objective values
close to the optimum. Including additional penalty terms in the
objective---small enough to not significantly affect the primary
objective---improves computation by encouraging binary variables to
take integral values in the relaxations, and also guides the solution
process towards particular solutions. For example, consider the
penalty terms (for a minimization problem)
\begin{equation}
\sum_{l \in \mc{L}} W^{y} (1-y_l) + \sum_{l \in \mc{L}} W^\ts{L}_l (1-\rho_l) + \sum_{g \in \mc{G}} W^\ts{G}_g (1-\zeta_l)
\end{equation}
where $W^y$, $W^\ts{L}_l$, $W^\ts{G}_g$ are weights to be chosen
appropriately. The first term penalizes $1-\cos\theta_l$, and hence
helps line loss. The second penalizes cuts to lines. For example,
setting $W^\ts{L}_l$ equal to the some small multiple of the
pre-islanding power flow through the line will penalize most heavily
disconnections of high-flow lines; in \cite{TBG+13} it was shown that
this leads more often to solutions that retain dynamic stability. The
third term penalizes the switching-off of generators. If
$W^\ts{G}_g=\epsilon P^{\ts{G}+}_g$ then units are given uniform
weighting. If, say, $W^\ts{G}_g =
\epsilon\bigl(P^{\ts{G}+}_g\bigr)^2$, then the disconnection of large
units is discouraged.

\subsection{Overall formulation}
The overall problem is to optimize the chosen islanding objective
(\eg~\eqref{eq:obj1}, \eqref{eq:obj2}, \eqref{eq:sumobj}, or
\eqref{eq:genmoveobj}), subject to
\begin{itemize}
\item sectioning constraints~\eqref{eq:IPsectioning};
\item line switching constraints~\eqref{eq:flow1};
\item power balance (\eqref{eq:linKCL1} and~\eqref{eq:linKCL2}) and flow~\eqref{eq:flow2} constraints;
\item the PWL approximation~\eqref{eq:xidef1};
\item generation limits~\eqref{eq:gencons};
\item line flow limits~\eqref{eq:linelim};
\item load shedding constraints~\eqref{eq:loadcons}.
\end{itemize}
\section{Computational results}
\label{sec:results}

\subsection{Islanding to minimize expected load loss}

A set of scenarios was built based on the $9$-, $14$-, $24$-, $30$-, $39$-, $57$-, $118$- and $300$-bus test systems from \textsc{Matpower}~\cite{MATPOWER}. For a network with $n^\ts{B}$ buses, $n^\ts{B}$
scenarios were generated by assigning in turn each single bus to
$\mc{B}^0$. No buses were included in $\mc{B}^1$ and no lines in
$\mc{L}^0$. For each scenario, the islanding solution was obtained by
solving the previously described MILP problem. The feasibility of an
islanding solution was checked by solving an AC optimal load shedding
(OLS) problem on the islanded network, which includes all AC power
balance, flow and operating constraints, but permits load shedding as
per~\eqref{eq:loadshedP} and~\eqref{eq:loadshedQ}.

Data for the islanding problems are described as follows. In the
objective function, $J^\ts{exp load}$, a value of $0.75$ is used for
the load loss penalty $\beta_d$. The generator coherency objective,
$J^\ts{coh}$, was not included initially. The penalties are
$W^\ts{G}_g = 0.01 P^{\ts{G}+}_g$, $W^y = 0.1$ and $W^\ts{L}_l =
0.0025 \sum_d P^\ts{D}_d$, so that the line-cut penalty is scaled by
the total load in the system. Our investigations show that these
penalties have a negligible effect ($0.2\%$) on the quality of the
solutions, but reduce computation time by an order of magnitude.  For
the PWL approximation for a line $l$, first the angle difference prior
to islanding, $\theta^*_l$, is determined from the base-case AC OPF
solution, and then $12$ pieces are used over $\pm \bigl( |\theta^*_l|
+ \unit{10}{\degree}\bigr)$.

Operating limits, including voltage and line limits, were obtained
from each network's data file~\cite{MATPOWER}. Generator real power
output limits ($P^{\ts{G}-}_g$ and $P^{\ts{G}+}_g$) were set, as
explained in Section~\ref{sec:opercons}, to allow a
  $2$-minute ramp change from the current output $P^{\ts{G}0}_g$,
  where ramp rates were available in the network data, or a $5\%$
  change where they were not. In either case, the output limits were
  limited by capacity limits. $P^{\ts{G}0}_g$ was obtained by solving
  an AC OPF on the intact network prior to islanding. Then in the
  islanding problem, the lower limit was raised by $5\%$ of
  $(P^{\ts{G}-}_g - P^{\ts{G}+}_g)$. The post-islanding AC OLS,
  however, was permitted to use the full range, $\bigl[ P^{\ts{G}-}_g,
    P^{\ts{G}+}_g \bigr]$. This avoids those solutions where an island
  is infeasible because of too much generated real power.

\subsubsection{AC-feasible islanding of 24-bus network}
\label{CompSwitchResults}

Returning to the example of Section~\ref{sec:infeas_example}, the PWL
AC islanding approach is applied to the problem of islanding bus
$6$. The islanding problem was solved both with and
  without the option (as part of the optimization) of switching the
  shunt reactor at bus $6$. The optimal solutions are shown in
  Fig.~\ref{fig:24bus}. Without shunt switching (PWL-AC-1), the cable
  $(6,10)$ is left intact and the final network topology is
  significantly different from before. With shunt switching permitted
  (PWL-AC-2), the cable is again switched, but fewer buses are
  islanded than for the DC solution.  The feasibility of each solution
  was checked by solving the AC OLS problem on the islanded network,
  and both PWL AC solutions satisfied all AC
  constraints. Tab.~\ref{tab:24bus} compares the DC, PWL-AC-1 and
  PWL-AC-2 solutions, using values obtained from both the MILP
  solutions and the post-islanding AC solutions. The PWL AC islanding
  solutions are close to the final AC OLS solutions. Note that the PWL
  AC solutions achieve AC feasibility at the cost of a lower expected
  load supply (hence higher expected load shed).

\begin{table}[b]
\centering\footnotesize
\caption{$24$-bus system: Comparison of solutions.}  
\begin{tabular}{lrrrr}
\toprule
Solution   				& DC 		& PWL-AC-1 	& PWL-AC-2  \\
\midrule
\multicolumn{4}{c}{MILP islanding solution}\\
$J^\ts{exp load}$ (MW)		& $2764.8$ & $2679.2$ & $2753.8$ \\
Generation (MW)				& $2850.0$ & $2892.7$ & $2844.1$ \\
Exp. load shed (MW)			&	$85.3$ &  $170.8$ & $96.2$	\\			
\midrule
\multicolumn{4}{c}{Post-islanding AC-OLS}\\
$J^\ts{exp load}$ (MW)		&	$\star$	& $2671.2$	& $2753.2$	\\
Generation 	(MW)			&	$\star$	& $2884.4$	& $2847.8$	\\				
Exp. load shed 	(MW)		&	$\star$	&  $178.9$	& $96.8$	\\
\bottomrule
\end{tabular}
\label{tab:24bus}
\end{table}

\subsubsection{Computation time}

The speed with which islanding decisions have to be made depends on
whether the decision is being made before a fault has occurred, as
part of contingency planning within secure OPF, or after, in which
case the time scale depends on the cause of the contingency.  Finding
solutions that are optimal, or to within a pre-specified percentage of
optimality, can take an unpredictable amount of time.  Hence,
especially in the latter case of reacting after a fault has occurred,
it is important to be able to produce good feasible solutions within
short time periods even if these are not necessarily optimal. To
illustrate how the quality of the solution depends on the solution
time, tests were run for a set of fixed times of between $5$ and $45$
seconds, returning the best found integer feasible
solution. Tab.~\ref{tab:gaps} summarizes these results for the $57$-,
$118$- and $300$-bus scenarios, quoting the average relative MIP gap
of returned solutions.  All the test cases with $39$ or fewer buses
solved to negligible \% gaps within 5 seconds, and are not shown.
Tab.~\ref{tab:gaps} also shows the average gaps between the returned
and best-known AC solutions for each scenario, where an AC solution
was obtained from a returned PWL islanding solution by solving the
AC-OLS on the islanded network. The mean error between the objectives
of the returned PWL-AC and AC solutions was less than
  $0.02\%$. For each network and scenario, the best-known
AC solution was the best from those found from the different
termination times, plus longer $1000$-second runs.  In the
  second and third sections of Tab.~\ref{tab:gaps}, the mean values
  are over all cases that were feasible within the time limit.  The
  platform was a $64$-bit Dual Intel Xeon processor and $128$~GiB RAM
  with up to 12 threads and using CPLEX $12.5$ as the MILP solver.
\begin{table}[t]
\centering\footnotesize
\caption{Solutions to islanding problems for different time limits.}\label{tab:gaps}%
\setlength{\tabcolsep}{8pt}
\begin{tabular}{lrrrrrr}
\toprule
Time (s) & $5$ & $10$ & $15$ & $20$ & $30$ & $45$ \\
\midrule
\multicolumn{7}{c}{Percentage with no islanding solution found within time} \\
$57$-bus & $1.7$ & $0.0$ & $0.0$ & $0.0$ &$0.0$ & $0.0$ \\
$118$-bus & $0.0$ & $0.0$ & $0.0$ & $0.0$  &$0.0$ & $0.0$ \\
$300$-bus & $17.7$ & $7.7$ & $3.0$ & $0.7$ & $0.0$ & $0.0$ \\
\midrule
\multicolumn{7}{c}{Mean \% between best MIP solution and the MIP bound} \\
$57$-bus & $0.04$ & $0.31$ & $0.21$ & $0.13$ & $0.13$ & $0.08$ \\
$118$-bus & $0.43$ & $0.13$ & $0.08$ & $0.05$ & $0.04$ & $0.03$ \\
$300$-bus & $0.08$ & $0.15$ & $0.11$ & $0.18$ & $0.14$ & $0.05$ \\
\midrule
\multicolumn{7}{c}{Mean \% between best AC solution found in time and best known} \\
$57$-bus & $0.19$ & $0.14$ & $0.12$ & $0.06$ &$0.06$ & $0.05$ \\
$118$-bus & $0.23$ & $0.09$ & $0.12$ & $0.04$ & $0.04$ & $0.04$ \\
$300$-bus & $0.03$ & $0.09$ & $0.07$ & $0.12$ & $0.11$ & $0.06$ \\
\bottomrule
\end{tabular}
\end{table}

The results show that good islanding solutions were found within
$\unit{30}{\second}$---and usually sooner---for all
networks. Moreover, the islanding topology usually changes little, or
not at all, between the solutions returned at $\unit{5}{\second}$ and
$\unit{45}{\second}$.

\subsubsection{AC feasibility}

Using the DC model $20\%$ of cases led to AC-infeasible
islands~\cite{TBG+13}, whereas none of the islands found using the PWL AC model
were infeasible.

\subsubsection{Promoting generator coherency}

The generator coherency objective, $J^\ts{coh}$, may be included for
the $24$-bus network example by taking second-order dynamic data taken
from~\cite{IEEERTS96}. For example, when $\mc{B}^0 = 3$, maximizing
just $J^\ts{exp load}$ leads to an optimal solution that places bus
$1$ in section $0$ along with bus $3$, and an expected load supply of
$2699$~MW. In doing this, the line between buses $1$ and $2$ is
switched, separating the large generator sets at these buses (which
would incur a cost of $J^\ts{coh} = 2.26$). However, when maximizing
the joint objective with $k = 100$, the optimal solution does not
include bus $1$ in section $0$, opting instead to leave the line
$(1,2)$ intact and placing just buses $3$ and $9$ in section $0$. With
$k = 100$, the expected load supply is slightly smaller ($2670$ MW),
but the strongly-coupled generators at buses $1$ and $2$ remain
connected ($J^\ts{coh} = 0.00$).
%

\subsection{Coherency-based islanding}

The coherency-based splitting approach was applied to the
$10$-machine, $39$-bus New England test network. Slow coherency
analysis, assuming second-order dynamics, shows that the machines may
be divided into two groups: those at buses $30$, $31$ and $39$
in one group, and then all others.

With $\mc{B}^0 = \{30,31,39\}$ and $\mc{B}^1 = \mc{B}^\ts{G} \setminus
\mc{B}^0$, the optimal solution splits the system as shown in
Tab.~\ref{tab:39bus}. Note that although buses $1$--$3$ and $5$--$9$
are included in the same section as $30$, $31$ and $39$, no
generators are present at these buses.  The objective was to minimize
the movement of generator real power outputs,
\ie~\eqref{eq:genmoveobj}. To achieve this split and leave the islands
balanced, the generator at bus $32$ has to lower its output from $671$
to $373$~MW, while $311$~MW is shed. It is worth stating that no other
solution exists that splits these two groups but
requires less total change in generator outputs.

\begin{table}[b]
\centering\footnotesize
\caption{Coherency-based islanding of
  $39$-bus network.}\label{tab:39bus}
\begin{tabular}{lrr}
\toprule
      & Section $0$ & Section $1$\\
\midrule
Buses & $1$--$3$, $5$--$9$, $30$, $31$, $39$ & $4$, $10$--$29$, $32$--$38$\\
Generation (MW) & $2007.18$ & $3992.29$ \\
Load supplied (MW) & $1997.89$ & $3945.37$ \\
Load shed (MW) & $297.21$ & $13.76$\\  
\bottomrule 
\end{tabular}
\end{table}

\section{Conclusions and future work}
\label{sec:conclusions}

An optimization-based framework for the intentional or
controlled islanding of power networks has been
presented. The approach is flexible with respect to the aims and
objectives of islanding, and finds islands that are balanced and
satisfy real and reactive power flow and operating constraints. It has
been shown that the inclusion of a piecewise linear model of AC power
flow allows AC-feasible islands to be found, where previously a
DC-based approach led to islands with out-of-bound voltages. The use
of objectives that promote generator coherency has been demonstrated.

Future work will investigate the wider practical aspects of the
approach by performing detailed simulations on representative networks
and blackout scenarios, considering transient and dynamic
performance. Current work is exploring the use of decomposition and
aggregation methods to improve the computational efficiency for larger
networks.

\bibliography{pwl_ac_refs}


\end{document}

%% file: trodd1.tikz
\pgfplotsset{
	cycle list={%
		{green!50!black, dashed, thick}, {red, thick}, {blue, dashdotted, thick}
	}%
}
\begin{tikzpicture}
\newcommand{\Vmax}{1.05}
\newcommand{\Vmin}{0.95}
\newcommand{\gL}{1}
\newcommand{\bL}{-5}
\newcommand{\bC}{1}
\newcommand{\Gii}{\gL}
\newcommand{\Gij}{(-1*\gL)}
\newcommand{\Gji}{(-1*\gL)}
\newcommand{\Gjj}{\gL}
\newcommand{\Bii}{(\bL + 0.5*\bC)}
\newcommand{\Bij}{(-1*\bL)}
\newcommand{\Bji}{(-1*\bL)}
\newcommand{\Bjj}{(\bL + 0.5*\bC)}
\newcommand{\ThetaMax}{40}
\newcommand{\ThetaMin}{-40}
\newcommand{\samples}{40}
\newcommand{\mywidth}{0.4\columnwidth}
\newcommand{\myheight}{4cm}

\begin{axis}[domain=\ThetaMin:\ThetaMax, samples=\samples, smooth, no markers, enlargelimits=false, legend pos=north west, ylabel={Real power (p.u.)}, width=\mywidth,height=\myheight, scale only axis=true, ylabel absolute, legend cell align=left, xlabel={$\theta_{ij}$ (deg)}, name=plot1, legend style={draw=none,fill=none}]
\newcommand{\pA}{%
	max(%
		\Gii*\Vmax*\Vmax + \Gij*\Vmax*\Vmin*cos(\x) + \Bij*\Vmax*\Vmin*sin(\x),
		\Gii*\Vmin*\Vmin + \Gij*\Vmin*\Vmin*cos(\x) + \Bij*\Vmin*\Vmin*sin(\x),
		\Gii*\Vmax*\Vmax + \Gij*\Vmax*\Vmax*cos(\x) + \Bij*\Vmax*\Vmax*sin(\x),
		\Gii*\Vmin*\Vmin + \Gij*\Vmin*\Vmax*cos(\x) + \Bij*\Vmin*\Vmax*sin(\x)%
		)
	}
\newcommand{\pB}{%
	min(%
		\Gii*\Vmax*\Vmax + \Gij*\Vmax*\Vmin*cos(\x) + \Bij*\Vmax*\Vmin*sin(\x),
		\Gii*\Vmin*\Vmin + \Gij*\Vmin*\Vmin*cos(\x) + \Bij*\Vmin*\Vmin*sin(\x),
		\Gii*\Vmax*\Vmax + \Gij*\Vmax*\Vmax*cos(\x) + \Bij*\Vmax*\Vmax*sin(\x),
		\Gii*\Vmin*\Vmin + \Gij*\Vmin*\Vmax*cos(\x) + \Bij*\Vmin*\Vmax*sin(\x)%
		)
	}
\newcommand{\perrA}{%
	max(
		(\Gii*\Vmax*\Vmax + \Gij*\Vmax*\Vmin*cos(\x) + \Bij*\Vmax*\Vmin*sin(\x)) - (\Gii*(2*\Vmax-1) + \Gij*(\Vmax + \Vmin + cos(\x) - 2) + \Bij*\x*pi/180),
		(\Gii*\Vmin*\Vmin + \Gij*\Vmin*\Vmin*cos(\x) + \Bij*\Vmin*\Vmin*sin(\x)) - (\Gii*(2*\Vmin-1) + \Gij*(\Vmin + \Vmin + cos(\x) - 2) + \Bij*\x*pi/180),
		(\Gii*\Vmax*\Vmax + \Gij*\Vmax*\Vmax*cos(\x) + \Bij*\Vmax*\Vmax*sin(\x)) - (\Gii*(2*\Vmax-1) + \Gij*(\Vmax + \Vmax + cos(\x) - 2) + \Bij*\x*pi/180),
		(\Gii*\Vmin*\Vmin + \Gij*\Vmin*\Vmax*cos(\x) + \Bij*\Vmin*\Vmax*sin(\x)) - (\Gii*(2*\Vmin-1) + \Gij*(\Vmin + \Vmax + cos(\x) - 2) + \Bij*\x*pi/180)%
		)
	}
\newcommand{\perrB}{%
	min(
		(\Gii*\Vmax*\Vmax + \Gij*\Vmax*\Vmin*cos(\x) + \Bij*\Vmax*\Vmin*sin(\x)) - (\Gii*(2*\Vmax-1) + \Gij*(\Vmax + \Vmin + cos(\x) - 2) + \Bij*\x*pi/180),
		(\Gii*\Vmin*\Vmin + \Gij*\Vmin*\Vmin*cos(\x) + \Bij*\Vmin*\Vmin*sin(\x)) - (\Gii*(2*\Vmin-1) + \Gij*(\Vmin + \Vmin + cos(\x) - 2) + \Bij*\x*pi/180),
		(\Gii*\Vmax*\Vmax + \Gij*\Vmax*\Vmax*cos(\x) + \Bij*\Vmax*\Vmax*sin(\x)) - (\Gii*(2*\Vmax-1) + \Gij*(\Vmax + \Vmax + cos(\x) - 2) + \Bij*\x*pi/180),
		(\Gii*\Vmin*\Vmin + \Gij*\Vmin*\Vmax*cos(\x) + \Bij*\Vmin*\Vmax*sin(\x)) - (\Gii*(2*\Vmin-1) + \Gij*(\Vmin + \Vmax + cos(\x) - 2) + \Bij*\x*pi/180)%
		)
	}
\newcommand{\pMerrA}{%
	max(
		(\Gii*\Vmax*\Vmax + \Gij*\Vmax*\Vmin*cos(\x) + \Bij*\Vmax*\Vmin*sin(\x)) - (\Gii*(2*\Vmax-1) + \Gij*(\Vmax + \Vmin + 1 - 2) + \Bij*\x*pi/180),
		(\Gii*\Vmin*\Vmin + \Gij*\Vmin*\Vmin*cos(\x) + \Bij*\Vmin*\Vmin*sin(\x)) - (\Gii*(2*\Vmin-1) + \Gij*(\Vmin + \Vmin + 1 - 2) + \Bij*\x*pi/180),
		(\Gii*\Vmax*\Vmax + \Gij*\Vmax*\Vmax*cos(\x) + \Bij*\Vmax*\Vmax*sin(\x)) - (\Gii*(2*\Vmax-1) + \Gij*(\Vmax + \Vmax + 1 - 2) + \Bij*\x*pi/180),
		(\Gii*\Vmin*\Vmin + \Gij*\Vmin*\Vmax*cos(\x) + \Bij*\Vmin*\Vmax*sin(\x)) - (\Gii*(2*\Vmin-1) + \Gij*(\Vmin + \Vmax + 1 - 2) + \Bij*\x*pi/180)%
		)
	}
\newcommand{\pMerrB}{%
	min(
		(\Gii*\Vmax*\Vmax + \Gij*\Vmax*\Vmin*cos(\x) + \Bij*\Vmax*\Vmin*sin(\x)) - (\Gii*(2*\Vmax-1) + \Gij*(\Vmax + \Vmin + 1 - 2) + \Bij*\x*pi/180),
		(\Gii*\Vmin*\Vmin + \Gij*\Vmin*\Vmin*cos(\x) + \Bij*\Vmin*\Vmin*sin(\x)) - (\Gii*(2*\Vmin-1) + \Gij*(\Vmin + \Vmin + 1 - 2) + \Bij*\x*pi/180),
		(\Gii*\Vmax*\Vmax + \Gij*\Vmax*\Vmax*cos(\x) + \Bij*\Vmax*\Vmax*sin(\x)) - (\Gii*(2*\Vmax-1) + \Gij*(\Vmax + \Vmax + 1 - 2) + \Bij*\x*pi/180),
		(\Gii*\Vmin*\Vmin + \Gij*\Vmin*\Vmax*cos(\x) + \Bij*\Vmin*\Vmax*sin(\x)) - (\Gii*(2*\Vmin-1) + \Gij*(\Vmin + \Vmax + 1 - 2) + \Bij*\x*pi/180)%
		)
	}

\addplot (\x,\pMerrA); 
\addplot (\x,\perrA); 
\addplot (\x,\pA);
\addplot (\x,\pMerrB); 
\addplot (\x,\perrB); 
\addplot (\x,\pB); 
\addlegendentry{Lin. error}
\addlegendentry{Lin.+cos error}
\addlegendentry{$p^{ij}_l$}
\end{axis} 
\begin{axis}[domain=\ThetaMin:\ThetaMax, samples=\samples, smooth, no markers, enlargelimits=false, xlabel={$\theta_{ij}$ (deg)}, ylabel={Reactive power (p.u.)}, width=\mywidth, height=\myheight, scale only axis=true, ylabel absolute, legend cell align=left, at=(plot1.right of north east), anchor=left of north west, legend style={draw=none,fill=none}]
\pgfplotsset{every axis legend/.append style={
	at={(0.95,0.98)},
	anchor=north east}}

\newcommand{\qA}{%
	max(
		-\Bii*\Vmax*\Vmax - \Bij*\Vmax*\Vmin*cos(\x) + \Gij*\Vmax*\Vmin*sin(\x),
		-\Bii*\Vmin*\Vmin - \Bij*\Vmin*\Vmin*cos(\x) + \Gij*\Vmin*\Vmin*sin(\x),
		-\Bii*\Vmax*\Vmax - \Bij*\Vmax*\Vmax*cos(\x) + \Gij*\Vmax*\Vmax*sin(\x),
		-\Bii*\Vmin*\Vmin - \Bij*\Vmin*\Vmax*cos(\x) + \Gij*\Vmin*\Vmax*sin(\x)%
		)
	}
\newcommand{\qB}{%
	min(
		-\Bii*\Vmax*\Vmax - \Bij*\Vmax*\Vmin*cos(\x) + \Gij*\Vmax*\Vmin*sin(\x),
		-\Bii*\Vmin*\Vmin - \Bij*\Vmin*\Vmin*cos(\x) + \Gij*\Vmin*\Vmin*sin(\x),
		-\Bii*\Vmax*\Vmax - \Bij*\Vmax*\Vmax*cos(\x) + \Gij*\Vmax*\Vmax*sin(\x),
		-\Bii*\Vmin*\Vmin - \Bij*\Vmin*\Vmax*cos(\x) + \Gij*\Vmin*\Vmax*sin(\x)%
		)
	}
\newcommand{\qerrA}{%
	max(
		(-\Bii*\Vmax*\Vmax - \Bij*\Vmax*\Vmin*cos(\x) + \Gij*\Vmax*\Vmin*sin(\x)) - (-\Bii*(2*\Vmax-1) - \Bij*(\Vmax + \Vmin + cos(\x) - 2) + \Gij*\x*pi/180),
		(-\Bii*\Vmin*\Vmin - \Bij*\Vmin*\Vmin*cos(\x) + \Gij*\Vmin*\Vmin*sin(\x)) - (-\Bii*(2*\Vmin-1) - \Bij*(\Vmin + \Vmin + cos(\x) - 2) + \Gij*\x*pi/180),
		(-\Bii*\Vmax*\Vmax - \Bij*\Vmax*\Vmax*cos(\x) + \Gij*\Vmax*\Vmax*sin(\x)) - (-\Bii*(2*\Vmax-1) - \Bij*(\Vmax + \Vmax + cos(\x) - 2) + \Gij*\x*pi/180),
		(-\Bii*\Vmin*\Vmin - \Bij*\Vmin*\Vmax*cos(\x) + \Gij*\Vmin*\Vmax*sin(\x)) - (-\Bii*(2*\Vmin-1) - \Bij*(\Vmin + \Vmax + cos(\x) - 2) + \Gij*\x*pi/180)%
		)
	}
\newcommand{\qerrB}{%
	min(
		(-\Bii*\Vmax*\Vmax - \Bij*\Vmax*\Vmin*cos(\x) + \Gij*\Vmax*\Vmin*sin(\x)) - (-\Bii*(2*\Vmax-1) - \Bij*(\Vmax + \Vmin + cos(\x) - 2) + \Gij*\x*pi/180),
		(-\Bii*\Vmin*\Vmin - \Bij*\Vmin*\Vmin*cos(\x) + \Gij*\Vmin*\Vmin*sin(\x)) - (-\Bii*(2*\Vmin-1) - \Bij*(\Vmin + \Vmin + cos(\x) - 2) + \Gij*\x*pi/180),
		(-\Bii*\Vmax*\Vmax - \Bij*\Vmax*\Vmax*cos(\x) + \Gij*\Vmax*\Vmax*sin(\x)) - (-\Bii*(2*\Vmax-1) - \Bij*(\Vmax + \Vmax + cos(\x) - 2) + \Gij*\x*pi/180),
		(-\Bii*\Vmin*\Vmin - \Bij*\Vmin*\Vmax*cos(\x) + \Gij*\Vmin*\Vmax*sin(\x)) - (-\Bii*(2*\Vmin-1) - \Bij*(\Vmin + \Vmax + cos(\x) - 2) + \Gij*\x*pi/180)%
		)
	}
\newcommand{\qMerrA}{%
	max(
		(-\Bii*\Vmax*\Vmax - \Bij*\Vmax*\Vmin*cos(\x) + \Gij*\Vmax*\Vmin*sin(\x)) - (-\Bii*(2*\Vmax-1) - \Bij*(\Vmax + \Vmin + 1 - 2) + \Gij*\x*pi/180),
		(-\Bii*\Vmin*\Vmin - \Bij*\Vmin*\Vmin*cos(\x) + \Gij*\Vmin*\Vmin*sin(\x)) - (-\Bii*(2*\Vmin-1) - \Bij*(\Vmin + \Vmin + 1 - 2) + \Gij*\x*pi/180),
		(-\Bii*\Vmax*\Vmax - \Bij*\Vmax*\Vmax*cos(\x) + \Gij*\Vmax*\Vmax*sin(\x)) - (-\Bii*(2*\Vmax-1) - \Bij*(\Vmax + \Vmax + 1 - 2) + \Gij*\x*pi/180),
		(-\Bii*\Vmin*\Vmin - \Bij*\Vmin*\Vmax*cos(\x) + \Gij*\Vmin*\Vmax*sin(\x)) - (-\Bii*(2*\Vmin-1) - \Bij*(\Vmin + \Vmax + 1 - 2) + \Gij*\x*pi/180)%
		)
	}
\newcommand{\qMerrB}{%
	min(
		(-\Bii*\Vmax*\Vmax - \Bij*\Vmax*\Vmin*cos(\x) + \Gij*\Vmax*\Vmin*sin(\x)) - (-\Bii*(2*\Vmax-1) - \Bij*(\Vmax + \Vmin + 1 - 2) + \Gij*\x*pi/180),
		(-\Bii*\Vmin*\Vmin - \Bij*\Vmin*\Vmin*cos(\x) + \Gij*\Vmin*\Vmin*sin(\x)) - (-\Bii*(2*\Vmin-1) - \Bij*(\Vmin + \Vmin + 1 - 2) + \Gij*\x*pi/180),
		(-\Bii*\Vmax*\Vmax - \Bij*\Vmax*\Vmax*cos(\x) + \Gij*\Vmax*\Vmax*sin(\x)) - (-\Bii*(2*\Vmax-1) - \Bij*(\Vmax + \Vmax + 1 - 2) + \Gij*\x*pi/180),
		(-\Bii*\Vmin*\Vmin - \Bij*\Vmin*\Vmax*cos(\x) + \Gij*\Vmin*\Vmax*sin(\x)) - (-\Bii*(2*\Vmin-1) - \Bij*(\Vmin + \Vmax + 1 - 2) + \Gij*\x*pi/180)%
		)
	}

\addplot (\x,\qMerrA); 
\addplot (\x,\qerrA); 
\addplot (\x,\qA);
\addplot (\x,\qMerrB); 
\addplot (\x,\qerrB); 
\addplot (\x,\qB); 
\addlegendentry{Lin. error}
\addlegendentry{Lin.+cos error}
\addlegendentry{$q^{ij}_l$}
\end{axis} 
\end{tikzpicture}

%% file: trodd2KM.tikz
%
%
%
%

\definecolor{mycolor1}{rgb}{0,0.75,0.75}

\begin{tikzpicture}

\begin{axis}[%
width=0.6\columnwidth,
height=6.25cm,
scale only axis=true,
xmin=65, xmax=100,
xlabel={Load (\% of peak)},
ymin=2.8, ymax=5.4,
ylabel={Generation cost (\$K/hour)},
legend pos = north west,
legend cell align=left,
legend style={draw=none}
]
\addplot [
color=blue,
dashed,
thick
]
table{
65 3.69397719201234
66 3.65048129917392
67 3.61537357606349
68 3.59158841884692
69 3.57216748222754
70 3.55671014729983
71 3.54598888605249
72 3.54112098757551
73 3.54384018861257
74 3.55713093727241
75 3.58709044854081
76 3.64278342488885
77 3.70362153873473
78 3.76515786718344
79 3.82739241005438
80 3.89032516719387
81 3.9539561384594
82 4.01828532370133
83 4.08331272279545
84 4.14903833560721
85 4.21546216198758
86 4.28288617675935
87 4.35126531704273
88 4.42037979273642
89 4.49023055043538
90 4.56081854411737
91 4.6321447357719
92 4.70421009456134
93 4.77701559777984
94 4.85055928416358
95 4.92483486621513
96 4.99984258118294
97 5.07558296795743
98 5.15205656542939
99 5.22926391300773
100 5.30720555065471
};
\addlegendentry{AC};

\addplot [
color=red,
solid,
thick
]
table{
65 3.88307554430931
66 3.82522546376586
67 3.77463250702908
68 3.73083128913523
69 3.6973448404382
70 3.67326613520534
71 3.65287508337068
72 3.63705301130449
73 3.62724237520372
74 3.62420421316381
75 3.63070226035108
76 3.65124765740662
77 3.69662448988688
78 3.75832625869371
79 3.82076088974823
80 3.88393658704745
81 3.94781809670871
82 4.01241296933905
83 4.07774023491215
84 4.14380044881894
85 4.21056742002228
86 4.27804230543445
87 4.34622807760885
88 4.41512946396225
89 4.48474672522244
90 4.55508012146122
91 4.62612991214702
92 4.6978963560495
93 4.77037971141764
94 4.84358189006729
95 4.91750113966932
96 4.99213829526603
97 5.06749362381276
98 5.14357280055783
99 5.22037905745669
100 5.2979046787583
};
\addlegendentry{SOS};

\addplot [
color=green!50!black,
dashdotted,
thick
]
table{
65 3.00657904731292
66 3.0596856101712
67 3.11346327112085
68 3.16800332457648
69 3.22335594552125
70 3.27938121033293
71 3.33609797078527
72 3.39344935924711
73 3.45145061898501
74 3.51045744305443
75 3.5700236432497
76 3.63022395069662
77 3.69123150614112
78 3.75294476405042
79 3.81537214878792
80 3.87854498363461
81 3.94244360397394
82 4.00704940356368
83 4.07236237636421
84 4.13838253093835
85 4.2051167057586
86 4.27257616847586
87 4.34076371963868
88 4.40967182162852
89 4.47930457505791
90 4.54965744239618
91 4.62073077634329
92 4.69252491857076
93 4.76504076332282
94 4.83829544421026
95 4.91227765082973
96 4.98698245154065
97 5.06241022704935
98 5.13861317851804
99 5.21557304175958
100 5.29325975151587
};
\addlegendentry{Relaxed};

\addplot [
color=violet,
densely dotted,
thick
]
table{
65 2.98402209547761
66 3.03616794595711
67 3.0889976615551
68 3.14251124227153
69 3.19670868810643
70 3.2515899990598
71 3.30715517513165
72 3.36340421632193
73 3.42033712263069
74 3.47795389405792
75 3.53625453060363
76 3.59523903226777
77 3.65490739905041
78 3.71525963095147
79 3.77629572797104
80 3.83801569010904
81 3.90041951736555
82 3.96350720974048
83 4.02727876723392
84 4.09173418984581
85 4.15687347757616
86 4.22269663042497
87 4.28920364839225
88 4.356394531478
89 4.42426927968223
90 4.49282789300489
91 4.56207037144603
92 4.63199671500564
93 4.70260692368372
94 4.77390099748026
95 4.84587893639526
96 4.91854074042873
97 4.99188640958069
98 5.06591594385105
99 5.14062934323996
100 5.2160266077473
};
\addlegendentry{DC};

\end{axis}
\end{tikzpicture}%